\documentclass[a4paper]{article}

\usepackage{natbib}
\usepackage{amssymb}
\usepackage{latexsym}
\usepackage{bussproofs}
\usepackage{url}

\usepackage{pxfonts}
\usepackage{microtype}%Improves spacing between letters and words, must be loaded after fonts

\usepackage{stackengine}
\stackMath

\newtheorem{definition}{Definition}

\newtheorem{theorem}{Theorem}
\newtheorem{corollary}{Corollary}

\begin{document}

\title{Normalisation and Subformula Property for a System of Intuitionistic Logic with General Introduction and Elimination Rules}
\author{Nils K\"urbis}
\date{}
\maketitle

\begin{center}
Published in \emph{Synthese} \\ 
\url{https://doi.org/10.1007/s11229-021-03418-8}\bigskip
\end{center}

\begin{abstract}
\noindent This paper studies a formalisation of intuitionistic logic by Negri and von Plato which has general introduction and elimination rules. The philosophical importance of the system is expounded. Definitions of `maximal formula', `segment' and `maximal segment' suitable to the system are formulated and corresponding reduction procedures for maximal formulas and permutative reduction procedures for maximal segments given. Alternatives to the main method used are also considered. It is shown that deductions in the system convert into normal form and that deductions in normal form have the subformula property.\bigskip 

\end{abstract}

\section{Introduction}
According to inferentialist semantics for the logical constants, the rules governing such an expression define its meaning, if those rules satisfy certain criteria. The view stems ultimately from Gentzen \citeyearpar[\S 5.13, 189]{gentzenuntersuchungen} and has received profound scrutiny and development in the hands of Dummett \citeyearpar{dummettphilbasisint, dummettLBM}, Prawitz \citeyearpar{prawitznaturaldeduction, prawitzideas, prawitzmeaningviaproofs,  prawitzdummett3} and Schroeder-Heister \citeyearpar{schroederheisterSEP}. There are intricate questions regarding the precise nature of these criteria and which logics may or may not fulfil them, for which I am here merely going to refer the reader to the quoted literature and references therein.\footnote{For a brief overview with a focus on Dummett and Prawitz, see \citep{kurbisPTSNM}.} The purpose of the present paper is not to contribute to their development or critique. In the light of its heritage and the vast literature ensuing its conception, an author is justified in taking for granted that inferentialist semantics for the logical constants is a topic worth exploring without needing to establish its foundations from scratch. I will here only rely on two aspects of the criteria that rules ought to satisfy that are generally agreed upon: (a) that there must be a certain balance, called, following Dummett, \emph{harmony}, between the grounds for deriving a complex formula by one of the introduction rules for its main operator and the consequences of using such a formula as the major premise of an elimination rule; (b) that deductions in a formalisation of logic in natural deduction ought to be subject to a \emph{normalisation theorem} which establishes that any formula that is the conclusion of an introduction rule and major premise of an elimination rule (for its main connective) may be removed from a deduction. (b) is generally regarded as a necessary condition for (a). A little more precisely, a deduction is in \emph{normal form} if (i) it contains no \emph{maximal formula}, that is a formula that is the conclusion of an introduction rule and the major premise of an elimination rule (for its main connective); and (ii) no \emph{maximal segment}, that is a sequence of formulas of the same shape arising from the applications of certain rules the last of which is major premise of an elimination rule.\footnote{\label{typetoken}I will, as best as readability allows, distinguish occurrences of formulas in deductions from a more abstract notion of formula that applies to formulas of the same shape or form, as it is customary to say. The latter could also be referred to as formula types, the former as their tokens. For brevity, by `formula' I often mean an occurrence of a formula, but I will be explicit about the distinction where this aids understanding. There are also schematic formulas and their instances, which may or may not be formulas of the same shape, by which the general statement of a rule of inference is distinguished from its application in a deduction: the former uses schematic formulas and specifies the common form of all its instances, the latter have formula occurrences as their premises and conclusions and are used in the construction of deductions. I will thus speak of rules as well as of their applications, but for the sake of brevity by `rule' I will often mean an application of a rule. This clarification and the ensuing greater precision in the use of terminology was added at the request of previous readers of this paper.} The normalisation theorem for a system establishes that any deduction in it can be transformed into a deduction in normal form. A deduction has the \emph{subformula property} iff any formula that occurs on it is a subformula of either an undischarged assumption or of the conclusion. If deductions in normal form have the subformula property, then for every deduction in the system, there is one with the subformula property. Details and precise definitions to follow in due course. The present paper proves a normalisation theorem for a system of intuitionistic logic by Negri and von Plato \citeyearpar[216f]{negriplatostructural} that captures aspect (a) of inferentialist semantics for the logical constants particularly neatly. It also satisfies aspect (b). Negri and von Plato do not prove normalisation for their system directly, but observe that it follows by translation into sequent calculus, a special case of a cut elimination theorem proved for a system of multiple conclusion sequent calculus by restriction to single conclusions, and translation back into natural deduction \citep[215]{negriplatostructural}. The subformula property follows therefrom. The contribution of the present paper is to provide a direct proof of normalisation for their system, which raises interesting issues and requires new techniques. Major consequences of normalisation, such as the subformula property and consistency, are also drawn. 

Negri and von Plato's system has some quite original features. The introduction rules for a connective $\ast$ are formulated in terms of the discharge of assumptions of the form $A\ast B$, and every rule of the system is one that allows the derivation of an arbitrary formula from side-deductions of that formula and some further premises, as is the case with disjunction elimination in Gentzen's system. The difference between introduction and elimination rules lies in whether a formula with the connective governed by the rules as main operator is discharged above a side-deduction required for an application of the rule or whether such a formula is a premise of the rule. 

Deductions in Gentzen's formalisation of intuitionistic logic in natural deduction can be brought into normal form and these deductions have the subformula property. This was shown by \citet[Ch. IV]{prawitznaturaldeduction}.\footnote{Von Plato has edited previously unpublished material of Gentzen's that shows that he had also proved these results for intuitionistic logic. See \citep{platogentzennormalisation} and \citep{platocellar}. A referee points out that Raggio published normalisation theorems for Gentzen's systems around the same time as Prawitz and remarked that `Gentzen has certainly had proofs' \citep[91]{raggiohauptsatz} of these results. Raggio's proof uses an interesting technique different from the one by now well known through the work of Prawitz. His normal form theorem is also a little different from Prawitz's. Raggio's method removes all segments that end in the major premise of an elimination rule at once, if \emph{all} their first formulas are concluded by an introduction rule or $\bot E$ (not restricted to atomic formulas).} Normalisation is a subtle process, and changing any rules in a system immediately raises the question whether both properties still hold for deductions after the modification. 

The present paper answers this question in the positive for Negri's and von Plato's formalisation of intuitionistic logic. To do so it is necessary to adjust the definition of `maximal formula', `segment' and `maximal segment' given by Prawitz, and accordingly to reformulate the reduction procedures to remove maximal formulas and maximal segments from deductions. This is new to the literature. The discussion in the main part of the paper is restricted to propositional logic. The conclusion considers suitable rules for the universal and existential quantifiers and equality, some of which are also new, and gives reduction procedures for maximal formulas with these expressions as main operators.

\section{General Elimination and Introduction Rules: Philosophical Considerations} 
Gentzen initially considers introduction and elimination rules for a primitive negation operator \citep[186]{gentzenuntersuchungen}, but then observes that, as $\neg A$ can be defined as $A\supset \bot$, they may be omitted \citep[189]{gentzenuntersuchungen}. The result is the system studied by \citet[Ch. IV]{prawitznaturaldeduction}: 

\begin{center}
\AxiomC{$A$} 
\AxiomC{$B$}
\LeftLabel{$\land I$: \ }
\BinaryInfC{$A\land  B$}
\DisplayProof\qquad\qquad
\AxiomC{$A\land  B$}
\LeftLabel{$\land E$: \ }
\UnaryInfC{$A$}
\DisplayProof\quad 
\AxiomC{$A\land  B$} 
\UnaryInfC{$B$}
\DisplayProof

\bigskip

\bottomAlignProof
\AxiomC{$[A]^i$}
\noLine
\UnaryInfC{$\Pi$}
\noLine
\UnaryInfC{$B$}
\LeftLabel{$\supset \!\! I$: \ }
\RightLabel{$_i$}
\UnaryInfC{$A\supset B$}
\DisplayProof \qquad \qquad
\bottomAlignProof
\AxiomC{$A\supset B$}
\AxiomC{$A$}
\LeftLabel{$\supset \!\! E$: \ }
\BinaryInfC{$B$}
\DisplayProof 

\bigskip

\bottomAlignProof
\AxiomC{$A$}
\LeftLabel{$\lor I$: \ }
\UnaryInfC{$A\lor B$}
\DisplayProof\quad
\bottomAlignProof
\AxiomC{$B$}
\UnaryInfC{$A\lor B$}
\DisplayProof\qquad\qquad
\bottomAlignProof
\AxiomC{$A\lor B$}
\AxiomC{$[A]^i$}
\noLine
\UnaryInfC{$\Pi$}
\noLine
\UnaryInfC{$C$}
\AxiomC{$[B]^i$}
\noLine
\UnaryInfC{$\Sigma$}
\noLine
\UnaryInfC{$C$}
\LeftLabel{$\lor E$: \ }
\RightLabel{$_i$}
\TrinaryInfC{$C$}
\DisplayProof

\bigskip

\AxiomC{$\bot$}
\LeftLabel{$\bot E$: \ } 
\UnaryInfC{$C$}
\DisplayProof
\end{center}

\noindent The elimination rules exhibit a certain discrepancy. The conclusions of applications of $\lor E$ and $\bot E$ can be any formula, whereas the conclusions of $\land E$ and $\supset \!\! E$ are subformulas of $A\land B$ and $A\supset  B$. 

The elimination rules for $\lor$ and $\bot$ do indeed belong in the same class: applying $\lor E$ requires two side-deductions, which also provide the conclusion of the rule, in which assumptions are discharged that are proper subformulas of the major premise $A\lor B$. $\bot$ has no proper subformulas, and correspondingly its elimination rule requires no side-deductions.\footnote{See \citep[Ch. 2.8]{kurbisproofandfalsity} for further discussion of the classification of these rules.} 

The discrepancy is eradicated in systems of natural deduction with \emph{general elimination rules} \citep[]{platogeneral}, which have alternative elimination rules for $\land$ and $\supset$ that are in the same class as those for $\lor$ and $\bot$: 
 
\begin{center}
\bottomAlignProof
\AxiomC{$A\land B$}
\AxiomC{$[A]^i, [B]^i$}
\noLine
\UnaryInfC{$\Pi$}
\noLine
\UnaryInfC{$C$}
\LeftLabel{$\land E$: \ }
\RightLabel{$_i$}
\BinaryInfC{$C$}
\DisplayProof\qquad 
\bottomAlignProof
\AxiomC{$A\supset B$}
\AxiomC{$A$}
\AxiomC{$[B]^i$}
\noLine
\UnaryInfC{$\Pi$}
\noLine
\UnaryInfC{$C$}
\LeftLabel{$\supset \!\! E$: \ }
\RightLabel{$_i$}
\TrinaryInfC{$C$}
\DisplayProof
\end{center}

\noindent From now on, the labels `$\land E$' and `$\supset\!\!  E$' refer to these two rules. 

Harmony between the introduction and elimination rules governing a logical constant $\ast$ consists in a certain balance between the grounds for deriving a formula with $\ast$ as main operator as specified by its introduction rules and the consequences that may be derived from such a formula as specified by its elimination rules. Harmony has two aspects. One is that the grounds for deriving $A\ast B$ as specified by $\ast I$ are balanced by its consequences as specified by $\ast E$; the other that the consequences of $A\ast B$ as specified by $\ast E$ are balanced by its grounds as specified by $\ast I$. 

Negri and von Plato explain that general elimination rules capture the thought that everything that follows from the grounds for deriving a formula $A\ast B$ follows from $A\ast B$. They name the following principle after a comparable one put forward by \citet[33f]{prawitznaturaldeduction}: 

\bigskip 

\noindent \textbf{Inversion Principle}: \emph{Whatever follows from the direct grounds for deriving a proposition must follow from that proposition}. \citep[6ff]{negriplatostructural}

\bigskip

\noindent In other words, all the consequences of the grounds of $A\ast B$ are consequences of $A\ast B$. Here the consequences of a proposition are consequences of that proposition together with any minor premises, as in the case of $\supset\!\!  E$ and $\lor E$. 

General elimination rules thus capture one aspect of harmony. What about the other aspect? We should expect it to be captured by a converse of this inversion principle: everything that follows from $A\ast B$ follows from the direct grounds for deriving $A\ast B$; all the consequences of $A\ast B$ are consequences of the grounds of $A\ast B$. By analogy, this should be captured by the introduction rules for $\ast$. But it is not immediate how Gentzen's rules might do so.

It is, however, immediate for Negri's and von Plato's \emph{general introduction rules} \citep[217]{negriplatostructural}, which have the rather original feature that, instead of introducing formulas with a connective $\ast$ as main operator as the conclusion of the rule, they permit discharge of assumptions of that form: 

\begin{center} 
\bottomAlignProof
\AxiomC{$A$}
\AxiomC{$B$}
\AxiomC{$[A\land B]^i$}
\noLine
\UnaryInfC{$\Pi$}
\noLine
\UnaryInfC{$C$}
\LeftLabel{$\land I$: \ }
\RightLabel{$_i$}
\TrinaryInfC{$C$}
\DisplayProof\qquad
\bottomAlignProof
\AxiomC{$[A]^i$}
\noLine
\UnaryInfC{$\Sigma$}
\noLine
\UnaryInfC{$B$}
\AxiomC{$[A\supset B]^i$}
\noLine
\UnaryInfC{$\Pi$}
\noLine
\UnaryInfC{$C$}
\LeftLabel{$\supset \!\! I$: \ }
\RightLabel{$_i$}
\BinaryInfC{$C$}
\DisplayProof\bigskip

\AxiomC{$A$}
\AxiomC{$[A\lor B]^i$}
\noLine
\UnaryInfC{$\Pi$}
\noLine
\UnaryInfC{$C$}
\LeftLabel{$\lor I$: \ }
\RightLabel{$_i$}
\BinaryInfC{$C$}
\DisplayProof\qquad
\AxiomC{$B$}
\AxiomC{$[A\lor B]^i$}
\noLine
\UnaryInfC{$\Pi$}
\noLine
\UnaryInfC{$C$}
\RightLabel{$_i$}
\BinaryInfC{$C$}
\DisplayProof
\end{center} 

\noindent From now on, the labels $\land I$, $\supset\!\!  I$ and $\lor I$ refer to these rules. They conform exactly to a converse of the inversion principle just quoted:\bigskip

\noindent \textbf{Converse Inversion Principle}: \emph{Whatever follows from a proposition must follow from the direct grounds for deriving that proposition}. 

\bigskip

\noindent Milne, another major figure in inferentialist semantics, concurs, albeit that he prefers classical over intuitionistic logic \citep{milneinversion}: he, too, proposes general introduction rules to capture harmony after having found the usual introduction rules stemming from Gentzen wanting in the light of considerations regarding the inversion principles.\footnote{\label{Milnefootnote}Milne formalised a system of classical logic with general introduction and elimination rules in which for every valid deduction, there is one with the subformula property \citep{milnesubformula}. This is rather unusual and a result of great importance. However, as my current purpose is not to adjudicate between classical and intuitionistic logic, I set Milne's system aside and focus on the logic favoured by Dummett and Prawitz. Milne's proof is model theoretic and not constructive. I establish his result proof-theoretically and constructively by means of normalisation in \citep{kurbisgenrulescl}, which can be read as a companion piece to the present paper. For the history of inversion principles, see \citep{moriconitesconi}.} 

General introduction and elimination rules, Negri and von Plato point out, exhibit a `perfect symmetry', captured in the following principle: 

\bigskip

\noindent \emph{General introduction rules state that if a formula $C$ follows from a formula $A$, then it already follows from the immediate grounds for $A$; general elimination rules state that if $C$ follows from the immediate grounds for $A$, then it already follows from $A$.} \citep[217]{negriplatostructural}

\bigskip

\noindent A system in which general elimination rules are paired with general introduction rules thus has a good claim on capturing harmony and its converse, a requirement Dummett calls \emph{stability} \citep[Ch. 13]{dummettLBM}. According to Dummett, stability between introduction and elimination rules is a necessary condition for those rules to define the meaning of the logical constant they govern. Negri and von Plato's system is thus of some philosophical importance: it arguably captures stability between introduction and elimination rules more accurately than Gentzen's and Prawitz's systems of intuitionistic logic and consequently the systems around which inferentialist semantics has centred. Critics of harmony, therefore, had better look at the present system. 

Negri and von Plato's system will be defined precisely in the next section.\footnote{It is an intuitionistic version of Parigot's classical system of free deduction when written in natural deduction style \citep{parigotfreededuction}. The latter is a system of multiple conclusion sequent calculus. The present system results if Parigot's pair of right conjunction rules are replaced by a single rule variant, the primitive negation rules are replaced by a suitable rule for $\bot$, the multiple conclusions are restricted to single conclusions in the most obvious way, and the result is transposed to the framework of natural deduction used here. This method is slightly simpler than the method by which Parigot constructs a more familiar version of natural deduction in sequent calculus style for classical logic from free deduction \citep[368]{parigotfreededuction}.} Gentzen's introduction and elimination rules for intuitionistic logic are easily derived from the general introduction and elimination rules (see below, p.\pageref{usualintrorules}). Conversely, the general introduction and elimination rules are easily derived from Gentzen's rules. Thus it is a formalisation of intuitionistic propositional logic.\bigskip

\noindent \emph{Remark 1. Primitive Negation.} The following is a pair of general introduction and elimination rules for a primitive negation operator: 

\begin{center}
\bottomAlignProof
\AxiomC{$[A]^i$}
\noLine
\UnaryInfC{$\Pi$}
\noLine
\UnaryInfC{$B$} 
\AxiomC{$[A]^i$}
\noLine
\UnaryInfC{$\Sigma$}
\noLine
\UnaryInfC{$\neg B$} 
\AxiomC{$[\neg A]^i$}
\noLine
\UnaryInfC{$\Xi$}
\noLine
\UnaryInfC{$C$}
\LeftLabel{$\neg I$: \ }
\RightLabel{$_i$}
\TrinaryInfC{$C$}
\DisplayProof\qquad
\bottomAlignProof
\AxiomC{$\neg A$}
\AxiomC{$A$}
\LeftLabel{$\neg E$: \ }
\BinaryInfC{$C$}
\DisplayProof
\end{center}

\noindent A familiar introduction rule for negation is derivable from $\neg I$ by letting $C$ be $\neg A$ and $\Xi$ empty, that is discharging $\neg A$ straight after assuming it. 

Notice that $\neg E$ does indeed have the form of a general elimination rule. This can be seen by returning to treating $\neg A$ as $A\supset \bot$ and replacing $\bot$ for $B$ in $\supset \!\! E$. As everything follows from $\bot$, a side-deduction showing that the conclusion of the application of the rule follows from $\bot$ is redundant. Analogously, treating $\neg$ as primitive again, $\neg E$ could be supplemented by a side-deduction of $C$ from two formulas $B$ and $\neg B$, which is redundant for the same reason (make use of the option of discharging only one assumption of the side-deduction and chose it to be identical to the conclusion: this derives the rule $\neg E$ above). These rules, however, are less satisfactory than those for the other connectives, as negation occurs in a premise of the introduction rule. This is often considered to be a shortcoming if rules are intended to define the meaning of a connective they govern, and for this reason it is preferable to define negation in terms of implication and \emph{falsum}.\footnote{For reflections on whether the meaning of negation is adequately defined by its usual rules, see \citep{kurbisnegation}. For an approach that justifies the negation of minimal logic, see \citep{kurbisminimal}.}\bigskip

\noindent \emph{Remark 2. Verum.} Dual to the \emph{falsum} constant $\bot$ is the \emph{verum} constant $\top$. Its general introduction rule allows its discharge at any moment in the deduction: 

\begin{center} 
\AxiomC{$[\top]^i$}
\noLine
\UnaryInfC{$\Pi$}
\noLine
\UnaryInfC{$C$}
\LeftLabel{$\top I$: \ }
\RightLabel{$_i$}
\UnaryInfC{$C$}
\DisplayProof
\end{center} 

\noindent $\top$ has no elimination rule.\bigskip  

\noindent \emph{Remark 3. Classical Logic.} Prawitz's formalisation of classical logic consists of Gentzen's rules for $\land$ and $\supset$ augmented by classical \emph{reductio ad absurdum}: 

\begin{center} 
\AxiomC{$[\neg P]^i$}
\noLine
\UnaryInfC{$\Pi$}
\noLine
\UnaryInfC{$\bot$}
\LeftLabel{$\bot_C$: \ }
\RightLabel{$_i$}
\UnaryInfC{$P$}
\DisplayProof 
\end{center}

\noindent where $P$ is an atomic formula.\bigskip 

\noindent The presence of this rule necessitates a restriction of the subformula property of deductions in normal form: allowance must be made for assumptions of the form $\neg P$ that are discharged by classical \emph{reductio ad absurdum} and formulas $\bot$ concluded from them \citep[42]{prawitznaturaldeduction}. Siders and von Plato prove a similar result for the full system of classical logic with general elimination rules \citep{platosiderclassical}.\footnote{In Milne's a system of classical logic, for every valid deduction, there is one with the \emph{unrestricted} subformula property. In the light of the necessity to restrict the subformula property in other formulations of classical logic, this is a remarkable result, but, for reasons given in footnote \ref{Milnefootnote}, I will not investigate it any further here.}

\section{Intuitionistic Propositional Logic} 
This section contains a more precise characterisation of the system \textbf{I} of intuitionistic logic with general introduction and elimination rules.\footnote{The definition of deductions in \textbf{I} follows the format used by Troelstra and Schwichtenberg \citeyearpar[Sec. 2.1.1]{troelstraschwichtenberg}.} The following section is the main section of this paper with the proof of the normalisation theorem and its corollaries. 

The definition of the language of \textbf{I} is standard. 

\begin{definition}[Connective, Atomic Formula, Degree of a Formula]\label{D1}
\normalfont $\bot$, $\supset$, $\land$ and $\lor$ are the \emph{connectives}. An \emph{atomic formula} is one that contains no connective. The \emph{degree} of a formula is the number of connectives occurring in it. 
\end{definition}

\noindent $\bot$, being a connective, is not an atomic formula, but a formula of degree 1.

Deductions in \textbf{I} have the familiar tree shape, with the (discharged or undischarged) assumptions at the top-most nodes or leaves and the conclusion at the bottom-most node or root. The conclusion of a deduction is said to depend on the undischarged assumptions of the deduction. Similar terminology is applied to subdeductions of deductions. 

Assumptions are assigned \emph{assumption classes}, (at most) one for each assumption, marked by a natural number, different numbers for different assumption classes. Formula occurrences of different types\footnote{See footnote \ref{typetoken}.} must belong to different assumption classes. Formula occurrences of the same type may, but do not have to, belong to the same assumption class. Discharge of assumptions is marked by a square bracket around the formula: $[A]^i$, $i$ being the assumption class to which $A$ belongs, with the same label also occurring at the application of the rule at which the assumption is discharged. Assumptions classes are chosen in such a way that if one assumption of an assumption class is discharged by an application of a rule, then it discharges all assumptions in that assumption class. Empty assumption classes are permitted: they are used in vacuous discharge, when a rule that allows for the discharge of assumptions is applied with no assumptions being discharged. 

Upper case Greek letters $\Sigma$, $\Pi$, $\Xi$, possibly with subscripts or superscripts, denote deductions. Often some of the assumptions and the conclusion of the deduction are mentioned explicitly at the top and bottom of $\Sigma$, $\Pi$, $\Xi$. Using the same designation more than once to denote subdeductions of a deduction means that these subdeductions are exact duplicates of each other apart from, possibly, the labels of the assumption classes: the deductions have the same structure and at every node formulas of the same type are premises and conclusions of applications of the same rules. 

\begin{definition}[Deduction in \textbf{I}]\label{D2}\ \\
\normalfont (i) The formula occurrence $A$ is a deduction in \textbf{I} of $A$ from the undischarged assumption $A$. 

\noindent (ii) If $\Sigma$, $\Pi$, $\Xi$ are deductions in \textbf{I}, then following are deductions of $C$ in \textbf{I} from the undischarged assumptions in $\Sigma$, $\Pi$, $\Xi$ apart from those in the assumption classes $i$ and $j$, which are discharged: 

\begin{center} 
\bottomAlignProof
\AxiomC{$\Sigma$}
\noLine
\UnaryInfC{$A$}
\AxiomC{$\Pi$}
\noLine
\UnaryInfC{$B$}
\AxiomC{$[A\land B]^i$}
\noLine
\UnaryInfC{$\Xi$}
\noLine
\UnaryInfC{$C$}
\RightLabel{$_{\land I \ i}$}
\TrinaryInfC{$C$}
\DisplayProof\qquad
\bottomAlignProof
\AxiomC{$\Sigma$}
\noLine
\UnaryInfC{$A\land B$}
\AxiomC{$[A]^i \ [B]^j$}
\noLine
\UnaryInfC{$\Pi$}
\noLine
\UnaryInfC{$C$}
\RightLabel{$_{\land E \ i, j}$}
\BinaryInfC{$C$}
\DisplayProof

\bigskip

\AxiomC{$[A]^i$}
\noLine
\UnaryInfC{$\Sigma$}
\noLine
\UnaryInfC{$B$}
\AxiomC{$[A\supset B]^j$}
\noLine
\UnaryInfC{$\Pi$}
\noLine
\UnaryInfC{$C$}
\RightLabel{$_{\supset I \ i, j}$}
\BinaryInfC{$C$}
\DisplayProof\qquad
\AxiomC{$\Pi$}
\noLine
\UnaryInfC{$A\supset B$}
\AxiomC{}
\noLine
\UnaryInfC{$\Sigma$}
\noLine
\UnaryInfC{$A$}
\AxiomC{$[B]^i$}
\noLine
\UnaryInfC{$\Xi$}
\noLine
\UnaryInfC{$C$}
\RightLabel{$_{\supset E \ i}$}
\TrinaryInfC{$C$}
\DisplayProof

\bigskip 

\AxiomC{}
\noLine
\UnaryInfC{$\Sigma$}
\noLine
\UnaryInfC{$A$}
\AxiomC{$[A\lor B]^i$}
\noLine
\UnaryInfC{$\Pi$}
\noLine
\UnaryInfC{$C$}
\RightLabel{$_{\lor I \ i}$}
\BinaryInfC{$C$}
\DisplayProof\quad
\AxiomC{}
\noLine
\UnaryInfC{$\Sigma$}
\noLine
\UnaryInfC{$B$}
\AxiomC{$[A\lor B]^i$}
\noLine
\UnaryInfC{$\Pi$}
\noLine
\UnaryInfC{$C$}
\RightLabel{$_{\lor I \ i}$}
\BinaryInfC{$C$}
\DisplayProof\qquad
\AxiomC{}
\noLine
\UnaryInfC{$\Sigma$}
\noLine
\UnaryInfC{$A\lor B$}
\AxiomC{$[A]^i$}
\noLine
\UnaryInfC{$\Pi$}
\noLine
\UnaryInfC{$C$}
\AxiomC{$[B]^j$}
\noLine
\UnaryInfC{$\Xi$}
\noLine
\UnaryInfC{$C$}
\RightLabel{$_{\lor E \ i, j}$}
\TrinaryInfC{$C$}
\DisplayProof

\bigskip

\AxiomC{$\Pi$}
\noLine
\UnaryInfC{$\bot$}
\RightLabel{$_{\bot E}$} 
\UnaryInfC{$C$}
\DisplayProof
\end{center}

\noindent (iii) Nothing else is a deduction in \textbf{I}.
\end{definition} 

\noindent We can suppress the label indicating the rule applied, but the labels indicating discharge must always be present. 

We may think of assumption classes as being assigned to formulas during the course of the construction of a deduction to mark the discharge of assumptions. Then the construction of deductions according to the definition leaves some assumptions without assumption classes in the completed deduction. We can assign them assumption classes afterwards.
\footnote{Troestra and Schwichtenberg write that assumptions `are supposed to be labeled by markers' \citep[36]{troelstraschwichtenberg} for assumption classes. We cannot decide at the outset which assumptions are discharged at which point during the construction of the deduction. But we can decide which ones are discharged by which application of a rule. I am grateful to a referee for pointing out errors in a previous attempt at defining deductions and consequence and suggestions for how to rectify them. Notice that if assumption classes are assigned to undischarged assumptions of completed deductions, then combining such deductions to form a new one requires deleting those labels for assumption classes. Similarly, applying the reduction procedures of the next section requires deleting square brackets enclosing discharged assumptions and their labels.} To record from which assumptions a conclusion has been derived, it then suffices to list the assumption classes to which the undischarged assumptions of the deduction belong. This will be a multiset. We can write $\Gamma\vdash_\mathbf{I} A$ if there is a deduction in \textbf{I} of (the formula occurrence) $A$ from (occurrences of) some of the formulas in $\Gamma$.\footnote{Structural rules for $\vdash$ follow: thinning by adding empty assumption classes of formulas, or splitting one assumption class into two, if it concerns formulas of the same type, contraction by relabelling two assumption classes with one of their labels.} 

The premise $A$ of $\supset \!\! E$ and $C$ in all three elimination rules are normally called the minor premises, but in the current system it is useful to have terminology that allows to distinguish them. 

\begin{definition}[Terminology for Premises and Discharged Assumptions]\label{majorminoretc}\ \\
\normalfont (i) In applications of the elimination rules, formula occurrences taking the places of $A\land B$, $A\supset B$, $A\lor B$ and $\bot$ to the very left above the line are the \emph{major premises}; formula occurrences taking the places of $C$ to their right are the \emph{arbitrary premises}, and a formula occurrence taking the place of $A$ inbetween in an application of $\supset \!\! E$ is the \emph{minor premise}. 

\noindent (ii) In applications of the introduction rules, formula occurrences taking the places of $A$ and $B$ to the very left above the line are the \emph{specific premises}, and those taking the place of $C$ to their right are the \emph{arbitrary premises}; formula occurrences taking the places of the discharged assumptions $A\supset B$, $A\lor B$ and $A\land B$ are \emph{the major assumptions discharged by} applications of the respective rules, and those taking the place of the discharged assumptions $A$ in $\supset\!\! I$ are the \emph{minor assumptions discharged by} applications of that rule.  
\end{definition}

\noindent Vacuous discharge happens when no assumption is discharged above an arbitrary premise or above the specific premise of $\supset \!\! I$. The latter is sometimes necessary, but the former is always superfluous: Instead of applying the rule, we might as well go on with the deduction straight from the arbitrary premise. In $\land E$, it is of course often necessary to make use of the option of discharging only one assumption. 

Applications of rules with vacuous discharge above arbitrary premises can be removed from deductions by what is often called \emph{simplification conversions}.\footnote{See, e.g., \citep[181]{troelstraschwichtenberg}.} As these procedures are obvious, I will give no details here. In the following, I will assume that any deduction is cleaned up so as to contain no vacuous discharge above arbitrary premises: vacuous discharge above arbitrary premises is banned. In particular, I will assume that this is done should vacuous discharge above an arbitrary premise arise as a result of the conversions of deductions that remove maximal formulas, to be given in the next section.\footnote{It is worth remarking that this cannot happen in Milne's formalisation of classical logic, where vacuous discharged may be banned altogether.} 

$\bot E$ can be restricted to atomic conclusions. The proof is by an induction over the degree of formulas and the following transformations, replacing the steps to the left of $\leadsto$ by those to its right: 

\begin{center}
\begin{tabular}{l c l}
\AxiomC{$\bot$} 
\UnaryInfC{$A\land B$}
\DisplayProof & \quad$\leadsto$\quad\quad & 
\AxiomC{$\bot$}
\UnaryInfC{$A$}
\AxiomC{$\bot$}
\UnaryInfC{$B$}
\AxiomC{$[A\land B]^1$}
\RightLabel{$_{\land I \ 1}$}
\TrinaryInfC{$A\land B$}
\DisplayProof\\
\\
\AxiomC{$\bot$}
\UnaryInfC{$A\lor B$}
\DisplayProof & \quad$\leadsto$\quad\quad & 
\AxiomC{$\bot$}
\UnaryInfC{$A$}
\AxiomC{[$A\lor B]^1$}
\RightLabel{$_{\lor I \ 1}$}
\BinaryInfC{$A\lor B$}
\DisplayProof\\
\\
\AxiomC{$\bot$}
\UnaryInfC{$A\supset B$}
\DisplayProof & \quad$\leadsto$\quad\quad & 
\AxiomC{$\bot$}
\UnaryInfC{$B$}
\AxiomC{$[A\supset B]^1$}
\RightLabel{$_{\supset I \ 1}$}
\BinaryInfC{$A\supset B$}
\DisplayProof
\end{tabular} 
\end{center} 

\noindent Obviously a step that concludes $\bot$ from $\bot$ by $\bot E$ is superfluous. From now on it is assumed that any application of $\bot E$ has an atomic conclusion.

\section{Normalisation for I} 
We begin by defining the notion of a maximal formula in a way that is suitable for the rules of the system \textbf{I}: 

\begin{definition}[Maximal Formula]\label{maximalformula}
\normalfont A \emph{maximal formula} with main operator $\ast$ in a deduction in \textbf{I} is an occurrence of a formula $A\ast B$ that is the major premise of an application of $\ast E$ and the major assumption discharged by an application of $\ast I$. 
\end{definition}

\noindent \emph{Reduction Procedures for Maximal Formulas} 

\noindent Maximal formulas are removed from deductions by applying the following \emph{reduction procedures for maximal formulas}, where $\Pi, \Sigma$ above $[A], [B]$ indicate that these deductions are used to conclude each formula occurrence in the assumption class to which $A, B$ belong (assumption class markers and square brackets are deleted). I will call the deduction to which a reduction procedure is applied \emph{the initial deduction} and the result of the conversion \emph{the reduced deduction}.\bigskip 

\noindent 1. The maximal formula has the form $A\land B$. Convert the deduction on the left into the deduction on the right: 

\begin{center}
\AxiomC{$\Sigma_1$}
\noLine
\UnaryInfC{$A$}
\AxiomC{$\Sigma_2$}
\noLine
\UnaryInfC{$B$}
\AxiomC{$[A\land B]^k$}
\AxiomC{$[A]^i \ [B]^j$}
\noLine
\UnaryInfC{$\Pi_1$}
\noLine
\UnaryInfC{$C$}
\RightLabel{$_{i, j}$}
\BinaryInfC{$C$}
\noLine
\UnaryInfC{$\Pi_2$}
\noLine
\UnaryInfC{$D$}
\RightLabel{$_k$}
\TrinaryInfC{$D$}
\DisplayProof\quad$\leadsto$\quad
\alwaysNoLine
\AxiomC{$\Sigma_1$}
\UnaryInfC{$A$}
\AxiomC{$\Sigma_2$}
\UnaryInfC{$B$}
\AxiomC{$\mathbin{\stackon[6pt]{[A]}{\Sigma_1}} \ \mathbin{\stackon[6pt]{[B]}{\Sigma_2}}$}
\UnaryInfC{$\Pi_1$}
\UnaryInfC{$C$}
\UnaryInfC{$\Pi_2$}
\UnaryInfC{$D$}
\RightLabel{$_k$}
\singleLine
\TrinaryInfC{$D$}
\DisplayProof
\end{center}

\noindent If assumption class $k$ contains only one formula (that is, the maximal formula removed by the procedure), then the final step by $\land I$ in the deduction to the right is omitted: in this case, the reduction procedure consists in replacing the deduction on the left only by the deduction that concludes $D$ by $\Pi_2$ from $\Sigma_1, \Sigma_2$ through $\Pi_1$ (that is, by the subdeduction concluding the arbitrary premise of $\land I$ on the right). The purpose of the final application of $\land I$ in the reduced deduction is to ensure that any other formulas in assumption class $k$ remain discharged after the application of the reduction procedure. If only the displayed maximal formula is in $k$, this purpose is not fulfilled and the application of $\land I$ introduces vacuous discharge; hence we omit it. Notice that applying the reduction procedure cannot introduce any new maximal formulas into the deduction. It can only introduce new maximal segments. More on this below, before the proof of the normalisation theorem.\footnote{A suggestion by a referee lead to an improvement in the description of the reduction procedure.}\bigskip

\newpage

\noindent 2. The maximal formula has the form $A\supset B$. Convert the deduction on the left into the deduction on the right: 

\begin{center} 
\AxiomC{$[A]^j$}
\noLine
\UnaryInfC{$\Sigma_1$}
\noLine
\UnaryInfC{$B$}
\AxiomC{$[A\supset B]^k$}
\AxiomC{$\Pi_1$}
\noLine
\UnaryInfC{$A$}
\AxiomC{$[B]^i$}
\noLine
\UnaryInfC{$\Pi_2$}
\noLine
\UnaryInfC{$C$}
\RightLabel{$_i$}
\TrinaryInfC{$C$}
\noLine
\UnaryInfC{$\Pi_3$}
\noLine
\UnaryInfC{$D$}
\RightLabel{$_{j, k}$}
\BinaryInfC{$D$}
\DisplayProof\quad$\leadsto$\quad
\alwaysNoLine
\AxiomC{$[A]^j$}
\UnaryInfC{$\Sigma_1$}
\UnaryInfC{$B$}
\AxiomC{$\Pi_1$}
\UnaryInfC{$[A]$}
\UnaryInfC{$\Sigma_1$}
\UnaryInfC{$[B]$}
\UnaryInfC{$\Pi_2$}
\UnaryInfC{$C$}
\UnaryInfC{$\Pi_3$}
\UnaryInfC{$D$}
\singleLine
\RightLabel{$_{j, k}$}
\BinaryInfC{$D$}
\DisplayProof
\end{center} 

\noindent As in the previous case, if assumption class $k$ contains only one formula (that is, the maximal formula removed by the procedure), then the final step by $\supset \!\! I$ in the deduction to the right is omitted. Furthermore, in case $\supset\!\! I$ was applied with vacuous discharge above its specific premise, the conversion may introduce applications of rules with vacuous discharge above arbitrary premises. This happens if an assumption in $\Pi_1$ is discharged above an arbitrary premise of a rule in $\Pi_3$. It is assumed that these are removed as part of the reduction procedure. Notice that applying this reduction procedure, too, cannot introduce any new maximal formulas into the deduction, and can only introduce new maximal segments.\bigskip

\noindent 3. The maximal formula has the form $A\lor B$. Convert the deduction on the left into the deduction on the right: 

\begin{center} 
\AxiomC{$\Sigma_1$}
\noLine
\UnaryInfC{$A$}
\AxiomC{$[A\lor B]^k$}
\AxiomC{$[A]^i$}
\noLine
\UnaryInfC{$\Pi_1$}
\noLine
\UnaryInfC{$C$}
\AxiomC{$[B]^j$}
\noLine
\UnaryInfC{$\Pi_2$}
\noLine
\UnaryInfC{$C$}
\RightLabel{$_{i, j}$}
\TrinaryInfC{$C$}
\noLine
\UnaryInfC{$\Pi_3$}
\noLine
\UnaryInfC{$D$}
\RightLabel{$_k$}
\BinaryInfC{$D$}
\DisplayProof\quad$\leadsto$\quad
\alwaysNoLine
\AxiomC{$\Sigma_1$}
\UnaryInfC{$A$}
\AxiomC{$\Sigma_1$}
\UnaryInfC{$[A]$}
\UnaryInfC{$\Pi_1$}
\UnaryInfC{$C$}
\UnaryInfC{$\Pi_3$}
\UnaryInfC{$D$}
\singleLine
\RightLabel{$_k$}
\BinaryInfC{$D$}
\DisplayProof
\end{center} 

\noindent As in the two previous case, if assumption class $k$ contains only one formula (that is, the maximal formula removed by the procedure), then the final step by $\lor I$ in the deduction to the right is omitted, and applying this reduction procedure, too, cannot introduce any new maximal formulas into the deduction and can only introduce new maximal segments. Similarly for the case where the specific premise of $\lor I$ is $B$ concluded by $\Sigma_2$.\bigskip

\noindent This completes the reduction procedures for maximal formulas. 

\bigskip 

\noindent \emph{Alternative Procedures}

\noindent It is worth mentioning some other ways of dealing with the fact that the application of $\ast I$ that gives rise to maximal formulas $A\ast B$ may discharge more formulas than the maximal formulas above the arbitrary premise $D$ in the deductions marked by $\Pi_1, \Pi_2, \Pi_3$. The following gives three alternatives, each of which avoids the final step by $\ast I$ to discharge those open assumptions in the reduced deduction.\bigskip

\noindent I. One alternative would be to add deductions of $A\ast B$ wherever there is such an assumption, using $\Sigma_1, \Sigma_2$ to conclude the specific premises of $\ast I$. These deductions also demonstrate how to derive the usual introduction rules for $\land, \lor, \supset$ from the general introduction rules:\label{usualintrorules}

\begin{center} 
\AxiomC{$\Sigma_1$}
\noLine
\UnaryInfC{$A$}
\AxiomC{$\Sigma_2$}
\noLine
\UnaryInfC{$B$}
\AxiomC{$[A\land B]^i$}
\RightLabel{$_i$}
\TrinaryInfC{$A\land B$}
\DisplayProof\qquad
\AxiomC{$[A]^i$}
\noLine
\UnaryInfC{$\Sigma_1$}
\noLine
\UnaryInfC{$B$}
\AxiomC{$[A\supset B]^j$}
\RightLabel{$_{i, j}$}
\BinaryInfC{$A\supset B$}
\DisplayProof\qquad
\AxiomC{$\Sigma_1$}
\noLine
\UnaryInfC{$A$}
\AxiomC{$[A\lor B]^i$}
\RightLabel{$_i$}
\BinaryInfC{$A\lor B$}
\DisplayProof
\end{center} 

\noindent Doing so only generates new maximal formulas in case $A\ast B$ was a maximal formula in the initial deduction and thus does not increase the number of maximal formulas in the reduced deduction.\bigskip

\noindent II. A second solution employs the fact that applications of introduction rules may be restricted to discharge only one occurrence of a formula. Suppose, for instance, one wanted to discharge $n$ formula occurrences of the type $A\lor B$ by an application of $\lor I$: 

\begin{prooftree}
\AxiomC{$\Sigma$}
\noLine
\UnaryInfC{$A$}
\AxiomC{$\underbrace{[A\lor B]^i, [A\lor B]^i \ldots [A\lor B]^i}$}
\noLine
\UnaryInfC{$\Pi$}
\noLine
\UnaryInfC{$C$}
\RightLabel{$_i$}
\BinaryInfC{$C$}
\end{prooftree}

\noindent Then instead of making this one application of $\lor I$, one can apply it $n$ times: 

\begin{prooftree}
\AxiomC{$\Sigma$}
\noLine
\UnaryInfC{$A$}
\AxiomC{$\Sigma$}
\noLine
\UnaryInfC{$A$}
\AxiomC{$\Sigma$}
\noLine
\UnaryInfC{$A$}
\AxiomC{$\underbrace{[A\lor B]^1, [A\lor B]^2 \ldots [A\lor B]^n}$}
\noLine
\UnaryInfC{$\Pi$}
\noLine
\UnaryInfC{$C$}
\RightLabel{$_1$}
\BinaryInfC{$C$}
\RightLabel{$_2$}
\BinaryInfC{$C$}
\noLine
\UnaryInfC{$\vdots$}
\noLine
\UnaryInfC{$C$}
\RightLabel{$_n$}
\BinaryInfC{$C$}
\end{prooftree} 

\noindent The cases for the other connectives are similar. There are two options for implementing this strategy: the restriction may be made either as part of the construction of deductions, or any deduction to be normalised is first transformed into one that satisfies the restriction before the reduction procedures are applied. Either option works, as the system with the unrestricted introduction rules and the system with their restricted versions are evidently equivalent. Obviously any application of a restricted introduction rule is also a correct application of the unrestricted version, and the converse holds in virtue of the following:\bigskip

\noindent \textbf{Proposition.} Any deduction can be transformed into one in which every application of a general introduction rule discharges exactly one major assumption.\bigskip

\noindent \emph{Proof.} By the ban on vacuous discharge above arbitrary premises, the transformations indicated above and an induction over a suitable measure of the complexity of deductions, e.g. the number of applications of introduction rules discharging multiple formula occurrences of highest degree in a deduction. Take such an application such that no other such application stands above it in the deduction. Applying the transformation reduces the measure.\bigskip

\noindent In the light of this proposition one could implement what may be called the \emph{unique discharge convention} on introduction rules: every application of an introduction rule for $\ast$ discharges exactly one formula occurrence of the form $A\ast B$. This has some advantages. The conclusion of an application of an introduction rule in Gentzen's system obviously occurs exactly once in a deduction, so if the unique discharge convention is upheld, there is a straightforward correspondence between deductions in Gentzen's system with the general elimination rule for $\supset$ and in the present system with general introduction rules.\footnote{It is for this reason that the unique discharge convention is appealed to in \citep{kurbisgenrulescl}: it permits an easy transposition of Milne's system into a more standard system of classical logic with the subformula property. The alternative reduction procedures for removing maximal formulas from deductions in \textbf{I} used here could also be adapted to Milne's system.} However, it also has disadvantages, as it is fair to say that upholding the unique discharge assumption destroys the most striking features of general introduction rules. Be that as it may, any sequence of applications of introduction rules as in the example above can be collapsed into one application, so one could, after maximal formulas have been removed from a deduction satisfying the unique discharge convention, also simplify it again in that respect, thereby restoring the characteristic and original features of general introduction rules.\bigskip

\noindent III. The third, and most interesting, alternative strategy is based on the observation that if an application of $\ast I$ gives rise to more than one maximal formula, then they may all be removed at once by a reduction procedure which simultaneously concludes formulas of the form $A\ast B$ discharged by $\ast I$ that are not maximal by the deductions given in the first alternative strategy.\bigskip 

\noindent 1. The maximal formula has the form $A\land B$. Convert the deduction on the left into the deduction on the right: 

\begin{center}
\AxiomC{$\Sigma_1$}
\noLine
\UnaryInfC{$A$}
\AxiomC{$\Sigma_2$}
\noLine
\UnaryInfC{$B$}
\AxiomC{$[A\land B]^k$}
\AxiomC{$[A]^i \ [B]^j$}
\noLine
\UnaryInfC{$\Pi_1$}
\noLine
\UnaryInfC{$C$}
\RightLabel{$_{i, j}$}
\BinaryInfC{$C$}
\noLine
\UnaryInfC{$\Pi_2$}
\noLine
\UnaryInfC{$D$}
\RightLabel{$_k$}
\TrinaryInfC{$D$}
\DisplayProof\qquad$\leadsto$\qquad
\alwaysNoLine
\AxiomC{$\mathbin{\stackon[6pt]{[A]}{\Sigma_1}} \ \mathbin{\stackon[6pt]{[B]}{\Sigma_2}}$}
\UnaryInfC{$\Pi_1^\ast$}
\UnaryInfC{$C$}
\UnaryInfC{$\Pi_2^\ast$}
\UnaryInfC{$D$}
\DisplayProof
\end{center}

\noindent where if there is only one formula in assumption class $k$, then $\Pi_1^\ast=\Pi_1$ and $\Pi_2^\ast=\Pi_2$, and if there is more than one formula, then $\Pi_1^\ast, \Pi_2^\ast$ are obtained from $\Pi_1, \Pi_2$ in the following way:\bigskip

\noindent (a) for formulas $A\land B$ in assumption class $k$ that are maximal: delete the application of $\land E$ that has the formula as major premise as in the pattern displayed above, by moving directly from the rule that concludes its arbitrary premise to the rule applied to its conclusion and concluding all assumptions $A$, $B$ discharged by this rule by $\Sigma_1, \Sigma_2$. 

\noindent (b) for formulas $A\land B$ in assumption class $k$ that are not maximal: conclude them by the derivation of the first strategy displayed on p.\pageref{usualintrorules} using $\Sigma_1, \Sigma_2$.\bigskip

\noindent Notice that, as in our `official' reduction procedure for maximal formulas of the form $A\land B$, applying the alternative procedure cannot introduce any new maximal formulas into the deduction. It may introduce new maximal segments, but the comments to be made in due course on this possibility in relation to the official procedure apply here, too. \bigskip

\noindent 2. The maximal formula has the form $A\supset B$. Convert the deduction on the left into the deduction on the right: 

\begin{center} 
\AxiomC{$[A]^j$}
\noLine
\UnaryInfC{$\Sigma_1$}
\noLine
\UnaryInfC{$B$}
\AxiomC{$[A\supset B]^k$}
\AxiomC{$\Pi_1$}
\noLine
\UnaryInfC{$A$}
\AxiomC{$[B]^i$}
\noLine
\UnaryInfC{$\Pi_2$}
\noLine
\UnaryInfC{$C$}
\RightLabel{$_i$}
\TrinaryInfC{$C$}
\noLine
\UnaryInfC{$\Pi_3$}
\noLine
\UnaryInfC{$D$}
\RightLabel{$_{j, k}$}
\BinaryInfC{$D$}
\DisplayProof\qquad$\leadsto$\qquad
\alwaysNoLine
\AxiomC{$\Pi_1^\ast$}
\UnaryInfC{$[A]$}
\UnaryInfC{$\Sigma_1$}
\UnaryInfC{$[B]$}
\UnaryInfC{$\Pi_2^\ast$}
\UnaryInfC{$C$}
\UnaryInfC{$\Pi_3^\ast$}
\UnaryInfC{$D$}
\DisplayProof
\end{center} 

\noindent where $\Pi_1^\ast, \Pi_2^\ast, \Pi_3^\ast$ are analogous to the previous case: if $k$ contains only one formula, they are identical to $\Pi_1, \Pi_2, \Pi_3$, otherwise they are obtained by deleting applications of $\supset \!\! I$ discharging maximal formulas in $k$, concluding assumptions $B$ becoming undischarged by $\Pi_1$, $\Sigma_1$ as in the pattern displayed, and concluding all others by the relevant deduction of the first strategy. Further comments on vacuous discharge and new maximal formulas and segments apply as usual.\bigskip

\noindent 3. The maximal formula has the form $A\lor B$. Convert the deduction on the left into the deduction on the right: 

\begin{center} 
\AxiomC{$\Sigma_1$}
\noLine
\UnaryInfC{$A$}
\AxiomC{$[A\lor B]^k$}
\AxiomC{$[A]^i$}
\noLine
\UnaryInfC{$\Pi_1$}
\noLine
\UnaryInfC{$C$}
\AxiomC{$[B]^j$}
\noLine
\UnaryInfC{$\Pi_2$}
\noLine
\UnaryInfC{$C$}
\RightLabel{$_{i, j}$}
\TrinaryInfC{$C$}
\noLine
\UnaryInfC{$\Pi_3$}
\noLine
\UnaryInfC{$D$}
\RightLabel{$_k$}
\BinaryInfC{$D$}
\DisplayProof\qquad$\leadsto$\qquad
\alwaysNoLine
\AxiomC{$\Sigma_1$}
\UnaryInfC{$[A]$}
\UnaryInfC{$\Pi_1^\ast$}
\UnaryInfC{$C$}
\UnaryInfC{$\Pi_3^\ast$}
\UnaryInfC{$D$}
\DisplayProof
\end{center} 

\noindent with $\Pi_1^\ast, \Pi_3^\ast$ constructed analogously to the previous cases. Further comments apply here, too. Similarly for the case where the premise of $\lor I$ is $B$ concluded by $\Sigma_2$.\bigskip

\noindent This completes the discussion of alternatives.\bigskip

\noindent Applications of general introduction and elimination rules require deductions of arbitrary premises $C$ which also provide the conclusion of the application of the rule. They form part of sequences of formula occurrences of the same shape\footnote{See footnote 1.}:

\begin{definition}[Segment]\label{segment}
\normalfont A \emph{segment} is a sequence of formula occurrences $C_1\ldots C_n$ of the same shape in a deduction such that $n>1$, for all $i<n$, $C_i$ is an arbitrary premise of an application of a rule and $C_{i+1}$ is its conclusion, and $C_n$ is not an arbitrary premise of an application of a rule. 
\end{definition} 
 
\noindent The \emph{length} of a segment is the number of formula occurrences of which it consists, its \emph{degree} the degree of any such formula. As $C_1\ldots C_n$ are all of the same shape, I will speak of the formula (as a type) constituting the segment.\bigskip

\noindent\emph{Observation.} A consequence of the ban on vacuous discharge above arbitrary premises is that the first formula of a segment is an arbitrary premise discharged by an introduction rule, the conclusion of which is the second formula of the segment. The major, minor and specific premises of rules are either assumptions or the last formulas of segments.

\begin{definition}[Maximal Segment]\label{maximalsegment} 
\normalfont A \emph{maximal segment} is a segment the last formula of which is the major premise of an elimination rule. 
\end{definition} 

\noindent Maximal segments are removed from deductions by \emph{permutative reduction procedures}. Of these there are 24 in total, as the major premises of $\supset\!\! E$, $\lor E$, $\land E$ and $\bot E$ can be derived by six rules (i.e. as the conclusions of the introduction and elimination rules for $\supset$, $\lor$ and $\land$). I will only give some of the cases for $\supset \!\! E$ and $\bot E$ as examples, the others being similar.\bigskip 

\noindent 1. The major premise of $\supset \!\! E$ is derived by $\lor I$. Convert the deduction on the left into the deduction on the right: 

\begin{center}
\def\defaultHypSeparation{\hskip .1in}
\AxiomC{$\Pi_1$}
\noLine
\UnaryInfC{$A$}
\AxiomC{$[A\lor B]^i$}
\noLine
\UnaryInfC{$\Pi_2$}
\noLine
\UnaryInfC{$C\supset D$}
\RightLabel{$_i$}
\BinaryInfC{$C\supset D$}
\AxiomC{$\Sigma_1$}
\noLine
\UnaryInfC{$C$}
\AxiomC{$[D]^j$}
\noLine
\UnaryInfC{$\Sigma_2$}
\noLine
\UnaryInfC{$E$}
\RightLabel{$_j$}
\TrinaryInfC{$E$}
\DisplayProof\quad$\leadsto$\quad
\def\defaultHypSeparation{\hskip .1in}
\AxiomC{$\Pi_1$}
\noLine
\UnaryInfC{$A$}
\AxiomC{$[A\lor B]^i$}
\noLine
\UnaryInfC{$\Pi_2$}
\noLine
\UnaryInfC{$C\supset D$}
\AxiomC{$\Sigma_1$}
\noLine
\UnaryInfC{$C$}
\AxiomC{$[D]^j$}
\noLine
\UnaryInfC{$\Sigma_2$}
\noLine
\UnaryInfC{$E$}
\RightLabel{$_j$}
\TrinaryInfC{$E$}
\RightLabel{$_i$}
\BinaryInfC{$E$}
\DisplayProof
\end{center} 

\noindent 2. The major premise of $\supset\!\!  E$ is derived by $\supset \!\! I$. Convert the deduction on the left into the deduction on the right: 

\begin{center}
\def\defaultHypSeparation{\hskip .1in}
\AxiomC{$[A]^i$}
\noLine
\UnaryInfC{$\Pi_1$}
\noLine
\UnaryInfC{$B$}
\AxiomC{$[A\supset B]^j$}
\noLine
\UnaryInfC{$\Pi_2$}
\noLine
\UnaryInfC{$C\supset D$}
\RightLabel{$_{i, j}$}
\BinaryInfC{$C\supset D$}
\AxiomC{$\Sigma_1$}
\noLine
\UnaryInfC{$C$}
\AxiomC{$[D]^k$}
\noLine
\UnaryInfC{$\Sigma_2$}
\noLine
\UnaryInfC{$E$}
\RightLabel{$_k$}
\TrinaryInfC{$E$}
\DisplayProof\quad$\leadsto$\quad
\def\defaultHypSeparation{\hskip .1in}
\AxiomC{$[A]^i$}
\noLine
\UnaryInfC{$\Pi_1$}
\noLine
\UnaryInfC{$B$}
\AxiomC{$[A\supset B]^j$}
\noLine
\UnaryInfC{$\Pi_2$}
\noLine
\UnaryInfC{$C\supset D$}
\AxiomC{$\Sigma_1$}
\noLine
\UnaryInfC{$C$}
\AxiomC{$[D]^k$}
\noLine
\UnaryInfC{$\Sigma_2$}
\noLine
\UnaryInfC{$E$}
\RightLabel{$_k$}
\TrinaryInfC{$E$}
\RightLabel{$_{i, j}$}
\BinaryInfC{$E$}
\DisplayProof
\end{center} 

\noindent 3. The major premise of $\supset\!\!  E$ is derived by $\land E$. Convert the deduction on the left into the deduction on the right: 

\begin{center}
\def\defaultHypSeparation{\hskip .1in}
\AxiomC{$\Pi_1$}
\noLine
\UnaryInfC{$A\land B$}
\AxiomC{$[A]^i \ [B]^j$}
\noLine
\UnaryInfC{$\Pi_2$}
\noLine
\UnaryInfC{$C\supset D$}
\RightLabel{$_{i, j}$}
\BinaryInfC{$C\supset D$}
\AxiomC{$\Sigma_1$}
\noLine
\UnaryInfC{$C$}
\AxiomC{$[D]^k$}
\noLine
\UnaryInfC{$\Sigma_2$}
\noLine
\UnaryInfC{$E$}
\RightLabel{$_k$}
\TrinaryInfC{$E$}
\DisplayProof\quad$\leadsto$\quad
\def\defaultHypSeparation{\hskip .1in}
\AxiomC{$\Pi_1$}
\noLine
\UnaryInfC{$A\land B$}
\AxiomC{$[A]^i \ [B]^j$}
\noLine
\UnaryInfC{$\Pi_2$}
\noLine
\UnaryInfC{$C\supset D$}
\AxiomC{$\Sigma_1$}
\noLine
\UnaryInfC{$C$}
\AxiomC{$[D]^k$}
\noLine
\UnaryInfC{$\Sigma_2$}
\noLine
\UnaryInfC{$E$}
\RightLabel{$_k$}
\TrinaryInfC{$E$}
\RightLabel{$_{i, j}$}
\BinaryInfC{$E$}
\DisplayProof
\end{center} 

\newpage

\noindent 4. The major premise of $\bot E$ has been derived by $\land I$. Convert the deduction on the left into the deduction on the right: 

\begin{center} 
\AxiomC{$\Pi_1$}
\noLine
\UnaryInfC{$A$}
\AxiomC{$\Pi_2$}
\noLine
\UnaryInfC{$B$}
\AxiomC{$[A\land B]^i$}
\noLine
\UnaryInfC{$\Sigma$}
\noLine
\UnaryInfC{$\bot$}
\RightLabel{$_i$}
\TrinaryInfC{$\bot$}
\UnaryInfC{$E$}
\DisplayProof\qquad $\leadsto$\qquad
\AxiomC{$\Pi_1$}
\noLine
\UnaryInfC{$A$}
\AxiomC{$\Pi_2$}
\noLine
\UnaryInfC{$B$}
\AxiomC{$[A\land B]^i$}
\noLine
\UnaryInfC{$\Sigma$}
\noLine
\UnaryInfC{$\bot$}
\UnaryInfC{$E$}
\RightLabel{$_i$}
\TrinaryInfC{$E$}
\DisplayProof
\end{center} 

\noindent 5. The major premise of $\bot E$ has been derived by $\lor E$. Convert the deduction on the left into the deduction on the right: 

\begin{center} 
\AxiomC{$\Pi$}
\noLine
\UnaryInfC{$A\lor B$}
\AxiomC{$[A]^i$}
\noLine
\UnaryInfC{$\Sigma_1$}
\noLine
\UnaryInfC{$\bot$}
\AxiomC{$[B]^j$}
\noLine
\UnaryInfC{$\Sigma_2$}
\noLine
\UnaryInfC{$\bot$}
\RightLabel{$_{i, j}$}
\TrinaryInfC{$\bot$}
\UnaryInfC{$E$}
\DisplayProof\qquad $\leadsto$ \qquad
\AxiomC{$\Pi$}
\noLine
\UnaryInfC{$A\lor B$}
\AxiomC{$[A]^i$}
\noLine
\UnaryInfC{$\Sigma_1$}
\noLine
\UnaryInfC{$\bot$}
\UnaryInfC{$E$}
\AxiomC{$[B]^j$}
\noLine
\UnaryInfC{$\Sigma_2$}
\noLine
\UnaryInfC{$\bot$}
\UnaryInfC{$E$}
\RightLabel{$_{i, j}$}
\TrinaryInfC{$E$}
\DisplayProof
\end{center} 

\noindent The other 19 permutative reduction procedures pose no further complications. 

\begin{definition}[Normal Form]\label{normalform}
\normalfont A deduction is in \emph{normal form} if it contains neither maximal formulas nor maximal segments. 
\end{definition}

\noindent Repeated application of a permutative reduction procedure reduces the length of a maximal segment by permuting applications of elimination rules upwards in the deduction. As observed earlier, the first formula of a segment can only be one discharged by an introduction rule, and so repeated application of a permutative reduction procedure turns a maximal segment into a maximal formula. At the top of every maximal segment, there stands a maximal formula, so to speak. 

\begin{definition}[Rank of Deductions]
\normalfont The \emph{rank} of a deduction $\Pi$ is the pair $\langle d, l\rangle$, where $d$ is the highest degree of a maximal formula or maximal segment in $\Pi$ or $0$ if there is none, and $l$ is the sum of the sum of the lengths of maximal segments of highest degree and the number of maximal formulas in $\Pi$. $\langle d, l\rangle < \langle d', l'\rangle$ iff either (i) $d<d'$ or (ii) $d=d'$ and $l<l'$.
\end{definition} 

\noindent Applying reduction procedures for maximal formulas cannot introduce new maximal formulas into the reduced deduction, but it may increase the lengths of maximal segments that were in the initial deduction.\footnote{The alternative reduction procedures may, incidentally, shorten maximal segments, namely if $C$ or $D$ form part of one.} In the case of maximal formulas of form $A\land B$, this can happen if $\Sigma_1$ concludes $A$ or $\Sigma_2$ concludes $B$ with an elimination rule and some formula occurrence in the assumption class to which the formulas discharged by $\land E$ belong is the major premise of an elimination rule in $\Pi_1$. In the case of maximal formulas of the form $A\supset B$, this can happen if $\Sigma_1$ concludes $B$ with an elimination rule and some formula occurrence in the assumption class to which the formulas discharged by $\supset \!\! E$ belong is the major premise of an elimination rule in $\Pi_2$, or if $\Pi_1$ concludes $A$ with an elimination rule and some formula occurrence in the assumption class to which the minor assumptions discharged by $\supset \!\! I$ belong is the major premise of an elimination rule in $\Sigma_1$. Similarly for maximal formulas of the form $A\lor B$. 

Any maximal segment that suffers an increase in length as a result of a reduction procedure is, however, of lower degree than the maximal formula removed, as the formulas that constitute the segment are subformulas of the latter. Hence applying a reduction procedure for maximal formulas cannot increase the rank of a deduction.  

Applying a permutative reduction procedures cannot introduce new maximal segments into the reduced deduction, but it may increase the lengths of maximal segments that were in the initial deduction. In examples 1.-3. above, this would happen if $E$ is part of a maximal segment.\footnote{It cannot happen in examples 4. and 5., as $\bot E$ is restricted to atomic conclusions.} To ensure all maximal segments are removed from a deduction, the permutative reduction procedures must be applied with a strategy. 

Say that a deduction that already is in normal form can be converted into itself. Then we have the following: 

\begin{theorem}\label{normalisation}
Any deduction in \textbf{I} can be converted into a deduction in normal form. 
\end{theorem} 

\noindent\emph{Proof.} The theorem follows by the considerations of the paragraphs immediately preceding the theorem and an induction over the rank of deductions. Applying reduction procedures for maximal formulas cannot increase the rank of a deduction, and as a maximal formula is removed, applying a reduction procedure to a maximal formula of highest degree decreases the rank of the deduction. Permutative reduction procedures must be applied so as to avoid an increase of the lengths of segments of highest degree. This can be achieved by applying one to a maximal segment of highest degree such that no maximal segment of highest degree stands above it in the deduction. This reduces the rank of the deduction.

\begin{corollary}
If $\Gamma\vdash_\mathbf{I} A$, then there is a deduction in normal form with an occurrence of $A$ as the conclusion and occurrences of some of the formulas in $\Gamma$ as the undischarged assumptions.
\end{corollary}

\noindent \emph{Proof.} Immediate from theorem \ref{normalisation}.\bigskip 

\noindent If there is a deduction of $C$ from assumptions $A_1 \ldots A_n$, then the deduction in normal form into which it converts may retain only some of these assumptions: applying the reduction procedures for maximal formulas of the form $A\supset B$ removes the deduction of the minor premises of $\supset \!\! E$, if $\supset \!\! I$ discharges vacuously above the specific premise.\footnote{In Milne's classical system, the deduction in normal form proceeds from the same assumptions.}

\begin{theorem}\label{majorassumptions}
If $\Pi$ is a deduction in normal form, then all major premises of elimination rules are (discharged or undischarged) assumptions of $\Pi$.
\end{theorem} 

\noindent\emph{Proof.} By the form of deductions in normal form, as a result of the permutative reduction procedures.

\begin{definition}[Branch]  
\normalfont A \emph{branch} in a deduction is a sequence of formula occurrences $\sigma_1\ldots \sigma_n$ such that $\sigma_1$ is an assumption of the deduction that is neither discharged by an elimination rule nor the major assumption discharged by an introduction rule, $\sigma_n$ is either the conclusion of the deduction or the minor premise of $\supset\!\! E$, and for each $n>i$: if $\sigma_i$ is the major premise of an elimination rule other than $\bot E$, $\sigma_{i+1}$ is an assumption discharged by it, and if it is the major premise of $\bot E$, $\sigma_{i+1}$ is the conclusion of the rule; if $\sigma_i$ is the specific premise of an introduction rule, $\sigma_{i+1}$ is a major assumption discharged by it; and if $\sigma_i$ is an arbitrary premise (of an introduction or an elimination rule rule), $\sigma_{i+1}$ is the conclusion of the rule. 
\end{definition} 

\noindent Branches begin with a formula occurrence that is either an undischarged assumption of the deduction or a minor assumption discharged by $\supset \!\! I$. Taking the formulas on a branch that form segments as units, we can also say that a branch consists of a sequence of formulas or segments. 

\begin{corollary} 
If any major premises of elimination rules are on a branch in a deduction in normal form, then they precede any major assumptions discharged by introduction rules that are on the branch.
\end{corollary} 

\noindent \emph{Proof.} By theorem \ref{majorassumptions}, the major premises of elimination rules that occur on a branch in a deduction in normal form are assumptions. Hence they are not the last formulas of any segments, and in particular they are not the last formulas of any segments beginning with discharged major assumptions of introduction rules.\bigskip

\noindent It is a consequence of theorem \ref{majorassumptions} that in a deduction in normal form the major premises of elimination rules do not form parts of segments. A branch in a deduction in normal form typically begins with a sequence of major premises of elimination rules, such that the conclusion of the last of them is either the first formula on a segment ending in the specific premise of an introduction rule (if it is $\bot E$) or the second formula of such a segment (in all other cases), and continues with a sequence of segments the first formulas of which are major assumptions discharged by introduction rules. The first half of the branch is called the \emph{E-part}, its second half the \emph{I-part}. Separating them is the \emph{minimal formula} or \emph{minimal segment}. It is the first formula or the first segment the last formula of which is the specific premise of an introduction rule. If the last application of an elimination rule is $\bot E$, there is a minimal formula and it is $\bot$. If the last application is any other elimination rule, there is a minimal segment. Either part may be empty: some branches in normal deductions consist of only an E-part, some of only an I-part, and in the case of a deduction that consists of a single formula $A$, both parts are empty and there is only a minimal formula. Inspection of the rules shows that all formulas on the E-part are subformulas of the previous one, and all formulas of the I-part are subformulas of the subsequent one. 

\begin{definition}[Order of Branches]
\normalfont A branch has order $0$ if its last formula is the conclusion of the deduction; it has order $n+1$ if its last formula is the minor premise of an application of $\supset \!\! E$ the major premise of which is on a branch of order $n$.
\end{definition}

\noindent A branch of order $0$ is also called a \emph{main branch} in the deduction. 

\begin{definition}[Subformula Property]
\normalfont A deduction $\Pi$ of a conclusion $C$ from the undischarged assumptions $A_1\ldots A_n$ has the \emph{subformula property} if every formula on the deduction is a subformula either of $C$ or of $A_1\ldots A_n$. 
\end{definition}

\noindent For brevity we may speak of a segment being the premise, conclusion or discharged assumption of a rule if its last or first formula is the premise, conclusion or discharged assumption of that rule.  

\begin{theorem}\label{subformula}
Deductions in normal form have the subformula property.
\end{theorem} 

\noindent \emph{Proof.} By inspection of the rules and an induction over the order of branches. Consider a branch of order 0. The branch begins with a (possibly empty) sequence of major premises of elimination rules, going from major premise to assumption discharged by the elimination rule, until it reaches a specific premise of an introduction rule, and then continues with segments discharged by introduction rules, until it reaches the conclusion of the deduction. All formulas on the latter part of the branch are subformulas of the conclusion of the deduction. All formulas on the former part of the branch are subformulas either of an assumption that remains undischarged in the deduction, in which case they are subformulas of a formula that is an undischarged assumption of the deduction, or they are subformulas of a formula discharged by $\supset\!\! I$, in which case they are subformulas of a subformula of the conclusion. A branch that ends in the minor premise of $\supset\!\! E$ ends in a formula that is a subformula of a branch of lower order, and hence the theorem holds by induction. 

\begin{corollary}
For any deduction in \textbf{I}, there is a deduction of the same conclusion from some of its undischarged assumptions with the subformula property. 
\end{corollary} 

\noindent \emph{Proof.} By theorems \ref{normalisation} and \ref{subformula}.\bigskip

\noindent Finally, let a \emph{proof} be a deduction of \textbf{I} that has no undischarged assumptions. Elimination rules do not discharge assumptions above their major premises. Hence if in a deduction in normal form there is a main branch that does not have an I-part, it is not a proof. Contraposing and applying theorem \ref{normalisation} establishes: 

\begin{corollary}\label{proofs}
If there is a proof of $A$ in \textbf{I}, then there is one that ends with an application of an introduction rule. 
\end{corollary}

\noindent The usual further corollaries follow. For instance, \textbf{I} has the disjunction property: if $\vdash_\mathbf{I} A\lor B$, then either $\vdash_\mathbf{I} A$ or $\vdash_\mathbf{I} B$. 

\begin{corollary}
\textbf{I} is consistent. 
\end{corollary}

\noindent\emph{Proof.} Suppose there is a proof of $\bot$ in \textbf{I}. Then by corollary \ref{proofs}, there is a proof of $\bot$ that ends with an application of an introduction rule. But $\bot$ has no introduction rule. Hence there is no proof of $\bot$.

\section{Conclusion}
Negri and von Plato only formulate general elimination rules for the quantifiers \citep[64]{negriplatostructural}, but not general introduction rules. They also do not give rules for equality. To close this paper, I will briefly consider the formalisation of a full system of intuitionistic predicate logic with equality with general introduction and elimination rules.\footnote{The rules for $\exists$ and $=$ of this section are also found in \citep{milneinversion}. In \citep{kurbisgenrulescl} it is shown that deductions in normal form in Milne's system of classical predicate logic with $\exists$, but not with $\forall$, satisfy the subformula property.}  

The language has two disjoint sets of variables, the \emph{parameters} $a, b, c\ldots$ playing the role of free variables, and the variables to be bound by the quantifiers $x, y, z\ldots$, which do not occur free in formulas. The terms of the language are built up from the parameters, constant symbols and function symbols in the usual way. An expressions that is like a formula or a term, but containing free variables instead of parameters, is often called a pseudo-formula or a pseudo-term. 

$A^x_t$ is the result of substituting all occurrences of the variable $x$ in $A$ by the term $t$.\footnote{If $t$ is a pseudo-term, it is assumed that none of its free variables gets bound by the substitution, i.e. that $t$ is free for $x$ in $A$. But notice that, if the result of the substitution is to be a formula rather than a pseudo-formula, as is the case in the use made of this notation in rules of inference, we need not consider this possibility: we only need to consider terms, not pseudo-terms. Analogously for the next kind of substitution.} $\Xi_t^a$ is the result of substituting the term $t$ for the parameter $a$ throughout deduction $\Xi$. 

The elimination rule for the existential quantifier already has the form of general elimination rules. The general elimination rule for the universal quantifier has the same form but with a different use of terms: 

\begin{center} 
\AxiomC{$\Sigma$}
\noLine
\UnaryInfC{$\exists x A$}
\AxiomC{$[A^x_a]^i$}
\noLine
\UnaryInfC{$\Pi$}
\noLine
\UnaryInfC{$C$}
\RightLabel{$_{\exists E \ i}$}
\BinaryInfC{$C$}
\DisplayProof\qquad \qquad
\AxiomC{$\Sigma$}
\noLine
\UnaryInfC{$\forall x A$}
\AxiomC{$[A^x_t]^i$}
\noLine
\UnaryInfC{$\Pi$}
\noLine
\UnaryInfC{$C$}
\RightLabel{$_{\forall E \ i}$}
\BinaryInfC{$C$}
\DisplayProof
\end{center} 

\noindent where in $\exists E$, the parameter $a$ does not occur in $\exists xA$, nor in $C$, nor in any formulas undischarged in $\Pi$ except those of the assumption class $[A^x_a]$. 

The following are general introduction rules for the quantifiers:

\begin{center} 
\AxiomC{$\Sigma$}
\noLine
\UnaryInfC{$A^x_t$}
\AxiomC{$[\exists xA]^i$}
\noLine
\UnaryInfC{$\Pi$}
\noLine
\UnaryInfC{$C$}
\RightLabel{$_{\exists I \ i}$}
\BinaryInfC{$C$} 
\DisplayProof\qquad \qquad
\AxiomC{$\Sigma$}
\noLine
\UnaryInfC{$A^x_a$}
\AxiomC{$[\forall x A]^i$}
\noLine
\UnaryInfC{$\Pi$}
\noLine
\UnaryInfC{$C$}
\RightLabel{$_{\forall I \ i}$}
\BinaryInfC{$C$}
\DisplayProof
\end{center} 

\noindent where in $\forall I$, the parameter $a$ does not occur in undischarged assumptions of $\Pi$. 

Deductions in the system of intuitionistic predicate logic are defined by adding clauses for these four rules to the inductive step of definition \ref{D2}. 

It is worth remarking that the rules for both quantifiers have the same form and differ only with respect to the occurrences of terms and parameters and consequently where restrictions on parameters are imposed. 

The major premise of an application of $\exists E, \forall E$ is, as before, its leftmost premise, the other being its arbitrary premise. Similarly for the specific and arbitrary premises of $\exists I, \forall I$. The major assumptions discharged by their applications are the formulas taking the places of $\exists xA$ and $\forall xA$, respectively. The definition of `maximal formula' is as before, as is that of `maximal segment', except that segments now of course also arise by applications of the rules for the quantifiers in the evident way. 

As we have an unlimited amount of parameters at our disposal, we may adopt the convention that every application of $\exists E$ and $\forall I$ has its own parameter, so that the parameter of an application of $\forall I$ only occurs in its specific premise and the formulas from which it is derived, and the parameter of an application of $\exists E$ occurs only in the formulas in the assumption class discharged by it and formulas derived from them. Consequently, the parameter occurs only above the application of the rule in a deduction, and any application of $\exists E$ or $\forall I$ below it has a different parameter. Call this the \emph{parameter convention}. 

Inspection of the reduction procedures for the propositional connectives shows that, if the parameter convention is upheld, then any correct application of $\exists E$ or $\forall I$ in the initial deduction remains correct in the reduced deduction. 

The following are the reduction procedures for maximal formulas of the form $\exists x A$ and $\forall xA$, continuing the numbering of those for propositional logic:\bigskip

\noindent 4. The maximal formula has the form $\exists xA$. Convert the deduction on the left into the deduction on the right: 

\begin{center} 
\AxiomC{$\Sigma$}
\noLine
\UnaryInfC{$A^x_t$}
\AxiomC{$[\exists x A]^k$}
\AxiomC{$[A^x_a]^j$}
\noLine
\UnaryInfC{$\Xi$}
\noLine
\UnaryInfC{$C$}
\RightLabel{$_j$}
\BinaryInfC{$C$}
\noLine
\UnaryInfC{$\Pi$} 
\noLine
\UnaryInfC{$D$}
\RightLabel{$_k$}
\BinaryInfC{$D$} 
\DisplayProof\qquad$\leadsto$\qquad
\AxiomC{$\Sigma$}
\noLine
\UnaryInfC{$A^x_t$}
\AxiomC{$\Sigma$}
\noLine
\UnaryInfC{$[A^x_t]$}
\noLine
\UnaryInfC{$\Xi ^a_t$}
\noLine
\UnaryInfC{$C$}
\noLine
\UnaryInfC{$\Pi$} 
\noLine
\UnaryInfC{$D$}
\RightLabel{$_k$}
\BinaryInfC{$D$} 
\DisplayProof
\end{center}

\noindent where the procedure is again as in the cases of the propositional connectives: if assumption class $k$ contains only one formula, the final step by $\exists I$ is omitted. By the parameter convention, in the initial deduction $a$ occurs only in $\Xi$, and hence after its replacement by $t$ it disappears altogether from the reduced deduction, which therefore is a correct deduction.\bigskip

\noindent 5. The maximal formula has the form $\forall xA$. Convert the deduction on the left into the deduction on the right: 

\begin{center} 
\AxiomC{$\Xi$}
\noLine
\UnaryInfC{$A^x_a$}
\AxiomC{$[\forall x A]^k$}
\AxiomC{$[A^x_t]^j$}
\noLine
\UnaryInfC{$\Sigma$}
\noLine
\UnaryInfC{$C$}
\RightLabel{$_j$}
\BinaryInfC{$C$}
\noLine
\UnaryInfC{$\Pi$} 
\noLine
\UnaryInfC{$D$}
\RightLabel{$_k$}
\BinaryInfC{$D$} 
\DisplayProof\qquad $\leadsto$\qquad
\AxiomC{$\Xi$}
\noLine
\UnaryInfC{$A^x_a$}
\AxiomC{$\Xi ^a_t$}
\noLine
\UnaryInfC{$[A^x_t]$}
\noLine
\UnaryInfC{$\Sigma$}
\noLine
\UnaryInfC{$C$}
\noLine
\UnaryInfC{$\Pi$} 
\noLine
\UnaryInfC{$D$}
\RightLabel{$_k$}
\BinaryInfC{$D$} 
\DisplayProof
\end{center} 

\noindent \noindent where the procedure is again as in the cases of the propositional connectives: if assumption class $k$ contains only one formula, the final step by $\forall E$ is omitted. By the parameter convention, in the initial deduction $a$ occurs only in $\Xi$, and, as it is no longer present in $\Xi_t^a$ upon replacement by $t$, this is the only place where it remains in the reduced deduction, which therefore is a correct deduction.\bigskip

\noindent The additional permutative reduction procedures for maximal segments pose no further problems, and I will not give them. 

The results of the previous section go through as before, if substitution instances of formulas of the form $\forall xA$ and $\exists x A$ are counted amongst their subformulas. Corollary \ref{proofs} for intuitionistic predicate logic is used to prove the existence property: if there is a proof of $\exists x A$, then, for some term $t$, there is a proof of $A_t^x$. 

Finally, what would general introduction and elimination rules for equality be? Equality raises a number of philosophical questions, not the least, in the present context, whether it is a logical constant the meaning of which may be defined by the rules of inference governing it. I will not try to answer this question here, but the second set of rules to be given, which effectively capture Leibniz' definition of equality, do, it seems to me, have a fair claim on satisfying the criteria of inferentialist semantics.\footnote{For an in depth discussion of equality in inferentialist semantics, a survey of existing proposals and a novel approach, see \citep{indrzejczakequality}.} 

The following is a general elimination rule for $=$:

\begin{prooftree}
\AxiomC{$\Pi$}
\noLine
\UnaryInfC{$t_1=t_2$}
\AxiomC{$\Sigma$}
\noLine
\UnaryInfC{$P_{t_1}^x$}
\AxiomC{$[P_{t_2}^x]^i$}
\noLine
\UnaryInfC{$\Xi$}
\noLine
\UnaryInfC{$C$}
\RightLabel{$_{=E \ i}$}
\TrinaryInfC{$C$}
\end{prooftree} 

\noindent where $P$ is atomic. The general case follows by induction. To exclude trivial applications of $=\!\!E$, we may require $t_1$ and $t_2$ to be different. 

One option of a general introduction rule for equality follows the pattern of the rule for $\top$ considered in Remark 2 of section 2. An assumption of the form $t=t$ may be discharged at any moment in a deduction: 

\begin{prooftree}
\AxiomC{$[t=t]^i$}
\noLine
\UnaryInfC{$\Pi$}
\noLine
\UnaryInfC{$C$}
\RightLabel{$_{=I \ i}$}
\UnaryInfC{$C$}
\end{prooftree} 

\noindent The ban on vacuous discharge prevents futile applications of this rule with no further effect than to deduce formulas from themselves. 

The major premise of $=\!\! E$ is $t_1=t_2$, its minor premise is $P_{t_1}^x$, the major assumption discharged by $=\!\!I$ is $t=t$, and in both rules $C$ is the arbitrary premise. 

Deductions in intuitionistic predicate logic with equality are defined by extending the inductive step of the definition of deductions in intuitionistic predicate logic by clauses for these two rules. 

The additional permutative reduction procedures for these rules follow the usual pattern. If the requirement of the difference of $t_1$ and $t_2$ in $=\!\!E$ is imposed, there are no maximal formulas of the form $t_1=t_2$. Otherwise the reduction procedure is straightforward and also follows a by now familiar pattern. Replace the deduction on the left by the deduction on the right, where the final application of $=\!\!I$ is omitted, if $k$ contains only one formula: 

\begin{center} 
\AxiomC{$[t=t]^k$}
\AxiomC{$\Sigma$}
\noLine
\UnaryInfC{$P_t^x$}
\AxiomC{$[P_t^x]^j$}
\noLine
\UnaryInfC{$\Xi$}
\noLine
\UnaryInfC{$C$}
\RightLabel{$_j$}
\TrinaryInfC{$C$}
\RightLabel{$_k$}
\UnaryInfC{$C$}
\DisplayProof\qquad$\leadsto$\qquad
\alwaysNoLine
\AxiomC{$\Sigma$}
\UnaryInfC{$P_t^x$}
\UnaryInfC{$\Xi$}
\UnaryInfC{$C$}
\singleLine
\RightLabel{$_k$}
\UnaryInfC{$C$}
\DisplayProof
\end{center} 

\noindent As $P$ is atomic, the reduction procedure cannot introduce new maximal formulas into the reduced deduction. 

A slightly more original option for a general introduction rule for equality results by modifying a rule proposed by Read to fit the present framework. Read observes that, if we add predicate parameters to the language, then $t_1=t_2$ may be inferred if there is a deduction of $F_{t_2}^x$ from $F_{t_1}^x$ in which $F$ is a predicate parameter not occurring in any undischarged assumptions except $F_{t_1}^x$, and conversely, a deduction of $F_{t_1}^x$ from $F_{t_2}^x$ in which $F$ is a predicate parameter not occurring in any undischarged assumptions except $F_{t_2}^x$ \citep[116]{readidentity}. This captures one half of Leibniz' definition of $t_1=t_2$ as $\forall F ( F_{t_1}^x\equiv F_{t_2}^x)$ in inferentialist terms. Cast into the form of general introduction rules, the rule becomes: 

\begin{prooftree}
\AxiomC{$[F_{t_2}^x]^i$}
\noLine
\UnaryInfC{$\Pi$}
\noLine
\UnaryInfC{$F_{t_1}^x$}
\AxiomC{$[F_{t_1}^x]^j$}
\noLine
\UnaryInfC{$\Sigma$}
\noLine
\UnaryInfC{$F_{t_2}^x$}
\AxiomC{$[t_1=t_2]^k$}
\noLine
\UnaryInfC{$\Xi$}
\noLine
\UnaryInfC{$C$}
\RightLabel{$_{=I' \ i, j, k}$} 
\TrinaryInfC{$C$}
\end{prooftree} 

\noindent where the predicate parameter $F$ does not occur in any formulas undischarged in $\Pi$ and $\Sigma$ except in those of the assumption classes $[F_{t_2}^x]$ and $[F_{t_1}^x]$.

The major assumption discharged by this rule is $t_1=t_2$. Suitable terminology for formulas in assumption classes $[F_{t_1}^x]$ and  $[F_{t_2}^x]$ would be the \emph{parametric assumptions discharged} by $=\!\!I'$. Its specific premises are the conclusions $F_{t_1}^x$ and $F_{t_2}^x$ of $\Pi$ and $\Sigma$. $C$ is the arbitrary premise. 

Considerations of harmony should then lead us to adding a second elimination rule for equality symmetrical to the first: 

\begin{prooftree}
\AxiomC{$\Pi$}
\noLine
\UnaryInfC{$t_1=t_2$}
\AxiomC{$\Sigma$}
\noLine
\UnaryInfC{$P_{t_2}^x$}
\AxiomC{$[P_{t_1}^x]^i$}
\noLine
\UnaryInfC{$\Xi$}
\noLine
\UnaryInfC{$C$}
\RightLabel{$_{=E^S \ i}$}
\TrinaryInfC{$C$}
\end{prooftree} 

\noindent where $P$ is atomic. Together $=\!\!E$ and $=\!\! E^S$ capture the other half of Leibniz' definition of equality in inferential terms. 

It would be possible to dispense with one of the deductions $\Pi$ and $\Sigma$ in $=\!\!I'$: given, say, $\Sigma$, replacing $Fx$ by $x=t_1$, which is possible if $F$ is a parameter satisfying the conditions of an application of $=\!\!I'$, gives a deduction of $t_2=t_1$ from $t_1=t_1$; the latter is provable by a single application of $=\!\!I'$ using $F_{t_1}^x$ as both premises and discharged assumptions; applying $=\!\!E$ to the thus concluded $t_1=t_2$ with $F_{t_2}^x$ as minor premise and $F_{t_1}^x$ as arbitrary premise, discharged assumption and conclusion yields the deduction of $F_{t_1}^x$ from $F_{t_2}^x$ required for an application of $=\!\!I'$. If this is done, harmony demands that the introduction rule for equality should be paired with only one elimination rule; dispensing with $\Pi$, this should be $=\!\!E$, $=\!\! E^S$ being derivable from the symmetry of equality, which in turn is derivable from $=\!\!E$ and the reflexivity of equality. For philosophical reasons, however, it may be preferable to leave $=\!\!I'$ as it is: the single deduction is sufficient only in the presence of $=\!\!E$, and so $=\!\!I'$ could not be said to define the meaning of $=$, while its elimination rule merely exploits this meaning as so defined according to the inversion principle.\footnote{Read's rule, in fact, has only one of the deductions. Read shows the redundancy of a second deduction by replacing $F$ by $\neg F$ in the first and contraposing. This move is not available in intuitionistic logic. A similar comment to the above applies: it would then not appear to be its introduction rule alone that defines the meaning of $=$, as an appeal to (classical) negation is required. For a discussion of the effect of various forms of or restrictions imposed on Leibniz' definition of equality, see \citep[Sec. 3]{indrzejczakequality}.}

If deductions in intuitionistic predicate logic with equality are defined by $=\!\!I'$ instead of $=\!\! I$ in the inductive step, then the reduction procedure for maximal formulas of the form $t_1=t_2$ is less trivial. We may assume that a corresponding version of the parameter convention is upheld for predicate parameters with respect to applications of $=\!\!I'$. Transform the deduction on the left into the deduction on the right, where $\Sigma_P^F$ is the result of substituting the predicate parameter $F$ by $P$ throughout $\Sigma$, and as always, the final step by $=\!\!I'$ is omitted if assumption class $k$ contains only one formula: 

\begin{center}
\def\defaultHypSeparation{\hskip .1in}
\AxiomC{$[F_{t_2}^x]^i$}
\noLine
\UnaryInfC{$\Pi$}
\noLine
\UnaryInfC{$F_{t_1}^x$}
\AxiomC{$[F_{t_1}^x]^j$}
\noLine
\UnaryInfC{$\Sigma$}
\noLine
\UnaryInfC{$F_{t_2}^x$}
\AxiomC{$[t_1=t_2]^k$}
\AxiomC{$\Xi_1$}
\noLine
\UnaryInfC{$P_{t_1}^x$}
\AxiomC{$[P_{t_2}^x]^l$}
\noLine
\UnaryInfC{$\Xi_2$}
\noLine
\UnaryInfC{$C$}
\RightLabel{$_l$}
\TrinaryInfC{$C$}
\noLine
\UnaryInfC{$\Xi_3$}
\noLine
\UnaryInfC{$D$}
\RightLabel{$_{i, j, k}$} 
\TrinaryInfC{$D$}
\DisplayProof\qquad$\leadsto$\qquad
\def\defaultHypSeparation{\hskip .1in}
\AxiomC{$[F_{t_2}^x]^i$}
\noLine
\UnaryInfC{$\Pi$}
\noLine
\UnaryInfC{$F_{t_1}^x$}
\AxiomC{$[F_{t_1}^x]^j$}
\noLine
\UnaryInfC{$\Sigma$}
\noLine
\UnaryInfC{$F_{t_2}^x$}
\alwaysNoLine
\AxiomC{$\Xi_1$}
\UnaryInfC{$P_{t_1}^x$}
\UnaryInfC{$\Sigma_P^F$}
\UnaryInfC{$P_{t_2}^x$}
\UnaryInfC{$\Xi_2$}
\UnaryInfC{$C$}
\UnaryInfC{$\Xi_3$}
\UnaryInfC{$D$}
\singleLine
\RightLabel{$_{i, j, k}$} 
\TrinaryInfC{$D$}
\DisplayProof
\end{center} 

\noindent By the parameter convention for predicate letters, $F$ occurs only in $\Pi$ and $\Sigma$, from which it disappears after replacement by $P$, so the reduced deduction is a correct deduction. $P$ being atomic, the reduction procedure does not introduce new maximal formulas into the deduction. In case the maximal formula arises from an application of $=\!\!E^S$, $P_{t_1}^x$ and $P_{t_2}^x$ are interchanged, and the reduction procedure replaces $F$ by $P$ in $\Pi$ instead of in $\Sigma$. 

Permutative reduction procedures pose once more no further problems.

\bigskip

\noindent \textbf{Acknowledgements.} I would like to thank the referees for \emph{Synthese} for their thoughtful engagement and suggestions for improvement. The research that lead to the final version of this paper was funded by the Alexander von Humboldt Foundation. 

\bigskip

\setlength{\bibsep}{0pt}
\bibliographystyle{chicago}
\bibliography{KurbisGenIntroRulesInt}
\end{document}